\documentclass{article}
\usepackage{amssymb}
\usepackage{geometry}
\usepackage{amsfonts}
\usepackage{indentfirst}
\usepackage{mathenv}
\usepackage{amsmath}
\usepackage[english,francais]{babel}
\usepackage[T1]{fontenc}
\usepackage{authblk}
\newtheorem{theorem}{Theorem}[section]
\newtheorem{lemma}[theorem]{Lemma}
\newtheorem{e-proposition}[theorem]{Proposition}

\newtheorem{e-definition}[theorem]{Definition\rm}
\newtheorem{remark}{\bf Remark\/}

\newtheorem{theoreme}{Th\'eor\`eme}[section]

\newtheorem{definition}[theoreme]{Definition\rm}

\newcommand{\cqfd}
{%
\mbox{}%
\nolinebreak%
\hfill%
\rule{2mm}{2mm}%
\medbreak%
\par%
}
\title{On the cost of fast controls for some families of dispersive or parabolic equations in one space dimension}
\author[1]{Pierre Lissy\footnote{lissy@ann.jussieu.fr}\thanks{Work supported by the ERC advanced grant 266907 (CPDENL) of the 7th Research Framework Programme (FP7)}}
\affil[1]{UPMC Univ Paris 06, UMR 7598, Laboratoire Jacques-Louis Lions, F-75005, Paris, France.}
\begin{document}
\maketitle
\selectlanguage{english}
\begin{abstract}
In this paper, we consider the cost of null controllability for a large class of linear equations of parabolic or dispersive type in one space dimension in small time. By extending the work of Tenenbaum and Tucsnak in \emph{New blow-up rates for fast controls of Schrödinger and heat equations}, we are able to give precise upper bounds on the time-dependance of  the cost of fast controls when the time of control $T$ tends to $0$. We also give a lower bound of the cost of fast controls for the same class of equations, which proves the optimality of the power of $T$ involved in the cost of the control. These general results are then applied to treat notably the case of linear KdV equations and fractional heat or Schrödinger equations.
\end{abstract}

\section{Introduction}
\subsection{Presentation of the problem}
This paper is devoted to studying fast boundary controls for some evolution equations of parabolic or dispersive type, with the spatial derivative not necessarily of second order.
 
Let $H$ be an Hilbert space (the state space) and $U$ be another Hilbert space (the control space). Let $A:\mathcal D(A)\rightarrow H$ be a self-adjoint operator with compact resolvent, the eigenvalues (which can be assumed to be different from $0$ without loss of generality) are called $(\lambda_k)_{k\geqslant 1}$, the eigenvector corresponding to the eigenvalue $\lambda_k$ is called $e_k$. We assume that $-A$ generates on $H$ a strongly continuous semigroup $S:t\mapsto S(t)=e^{-tA}$. The Hilbert space $D(A^*)'(=D(A)')$ is from now on equipped with the norm
$$||x||_{D(A)'}^2=\sum \frac{<x,e_k>_H^2}{\lambda_k^2} .$$

We call $B\in \mathcal L_c(U,\mathcal D(A)')$ an admissible control operator for this semigroup, i.e. such that there exists some time $T_0>0$, there exists some constant $C>0$ such that for all $z\in \mathcal D(A)$, one has $$\int_0^{T_0}||B^*S(t)^*z||^2_U\leqslant C||z||^2_H.$$
We recall that if $B$ is admissible, then necessarily the previous inequality holds at every time, that is to say for all time $T>0$, there exists some constant $C(T)>0$ such that for all $z\in \mathcal D(A)$, one has $$\int_0^T||B^*S(t)^*z||^2_U\leqslant C(T)||z||^2_H.$$

From now on, we consider control systems of the following form:

\begin{gather}\label{heat}y_t+Ay=Bu\end{gather} or 
\begin{gather}\label{sch}y_t+iAy=Bu,\end{gather}
where $A$ will always be supposed to be positive in the parabolic case (i.e. for Equation \eqref{heat}).
Then, it is well-known (see for example \cite[Chapter 2, Section 2.3]{MR2302744}, the operators $-A$ or $-iA$ generates a  strongly continuous semigroup under the hypothesis given before thanks due to the Lummer-Phillips or Stone theorems) that if $u\in L^2((0,T),U)$, System \eqref{heat} or \eqref{sch} with initial condition $y^0\in H$ has a unique solution satisfying $y\in \mathcal C^0([0,T],H)$. Moreover, if the system is null controllable at some time $T_0$ (i.e. for all $y^0\in H$, there exists some control $u\in L^2((0,T_0),U)$ such that $y(T_0,\cdot)\equiv 0$), then there exists a unique optimal (for the $L^2((0,T_0),U)$-norm) control $u_{opt}\in L^2((0,T_0),U)$, the map $y^0\mapsto u_{opt}$ is then linear continuous (see for example \cite[Chapter 2, Section 2.3]{MR2302744}). The norm of this operator is called the optimal null control cost at time $T_0$ (or in a more concise form the cost of the control) and denoted $C_{T_0}$, which is also the smallest constant $C>0$ such that for all $y^0\in H$, there exists some control $u$ driving $y^0$ to $0$ at time $T_0$ with $$||u||_{L^2((0,T_0),U)}\leqslant C||y^0||_{H}.$$

Concerning \eqref{sch}, it can be shown (see for example \cite[Chapter 2, Section 2.3, Theorem 2.41]{MR2302744}) that this system is null controllable if and only if it is exactly controllable; moreover, in this case, it is easy to prove that the cost of exact controllability has the same behavior in small time as the cost of null controllability; hence, even for conservative systems, we will only be interested in null controllability.

Our goal in this work is to estimate precisely the cost of the control $C_T$ when the time $T\rightarrow 0$ for some families of operators $A$ which are null controllable in arbitrary small time, and notably to find lower and upper bounds on $C_T$. Understanding the behavior of fast controls is of great interest in itself but it may also be applied to study the uniform controllability  of transport-diffusion in the vanishing viscosity limit as explained in \cite{MR2956149}. (the strategy described in \cite{MR2956149} might probably be extended to the study of other problems of uniform controllability, for example in zero dispersion limit or in zero diffusion-dispersion limit as in \cite{MR2463799} or \cite{MR2571687})
It is obvious that $C_T$ must tend to $\infty$ when $T\rightarrow 0$. 

\subsection{State of the art}
In all what follows, $f\lesssim g$ (with $f$ and $g$ some complex valued functions depending on some variable $x$ in some set $\mathcal S$)  means that there exists some constant $C>0$ (possibly depending on other parameters) such that for all $x\in\mathcal S$, one has
$|f(x)|\leqslant C|g(x)|$, (such a $C$ is called an implicit constant in the inequality $f\lesssim g$), and $f\simeq g$  means that we have both $f\lesssim g$ and $g\lesssim f$. Sometimes, when it is needed, we might detail the dependance of the implicit constant with respect to some parameters. We also might write $g\gtrsim f$ when $f\lesssim g$. The set $\mathcal S$ will not be explicitly given, it will in general correspond to all the variables appearing explicitly in the inequality.

As far as we know, results concerning the cost of fast boundary controls have been obtained essentially in the case of heat and Schrödinger equations. It is known for a long time that for the one-dimensional heat equation posed on a time-space cylinder $(0,T)\times(0,L)$ with boundary control on one side, the time-dependence of the cost of the boundary control is $\simeq e^{\frac{\alpha^+}{T}}$ for some constant $\alpha>0$ (see \cite{MR805273} and \cite{MR743923}), where the notation $\alpha^+$ means that we simultaneously have that the cost of the control is $\gtrsim e^{\frac{\alpha}{T}}$ and $\lesssim e^{\frac{K}{T}}$ for all $K>\alpha$ as close as $\alpha$ as we want (the implicit constant in front of the exponential might possibly explode when we get closer to $\alpha$). The constant $\alpha$ verifies 
 $$L^2/4\leqslant \alpha\leqslant 3L^2/4.$$ The upper bound is obtained in \cite{TT} and the lower bound in \cite{2004-Miller} (it is the best bounds obtained until now). For the Schrödinger equation posed on a time-space cylinder $(0,T)\times(0,L)$ with boundary control on one side, one also has that
 the dependence in time of the cost of the boundary control is under the form $\simeq e^{\frac{\tilde\alpha^+}{T}}$ for some constant $\tilde\alpha>0$. The constant $\tilde\alpha$ verifies 
 $$L^2/4\leqslant \tilde\alpha\leqslant 3L^2/2.$$ The upper bound is obtained in \cite{TT} and the lower bound in \cite{MR2062431} (it is the bounds obtained until now). In both cases, it is conjectured that the lower bound is optimal, i.e. that one can choose
$$\alpha=\tilde\alpha=L^2/4.$$ 

Our goal is to extend this results to other first-order time evolution equations with spatial operators that are self-adjoint or skew-adjoint with eigenvalues $\lambda_k$ or $i\lambda_k$ that do not necessarily behave has $k^2$ or $ik^2$, for example linear KdV equations, anomalous diffusion equations or fractional Schrödinger equations. Our main tool is the \emph{moment method} which was introduced in \cite{71FR} for the study of heat-like equations in one space dimension (and more generally for parabolic systems with eigenvalues having a behavior as in equation \eqref{puiss} for some $\alpha>1$) and used successfully many times notably to prove the controllability or uniform controllability of parabolic systems or equations or to study the behavior of the cost of the control (see for example \cite{05aa}, \cite{TT}, \cite{2010-Glass-JFA} or \cite{MR2851682}). We first prove some general theorems about the cost of controls for operators $A$ having eigenvalues which behave asymptotically as $k^{\alpha}$ for some $\alpha\geqslant 2$, and give precise upper bounds concerning $C_T$. Concerning lower bounds, we also prove that $\limsup_{T\rightarrow 0} T^{1/(\alpha-1)}\ln(C_T)>0$ as soon as $\alpha>1$. These main theorems are then applied to some families of equations, as described further. However, since our work is mainly an extension of \cite{TT}, we are not going to improve any existing upper bounds in the case of heat or Schrödinger equations.

Concerning linear dispersive equations of KdV type, the controllability has been widely studied with different boundary conditions and different boundary controls (see, in particular, \cite{MR1214759}, \cite{MR1360229}, \cite{97Ro}, \cite{MR2084328}, \cite{MR2571687}, \cite{MR2463799} or \cite{MR2724598}), in general in order to prove a result of controllability for the corresponding nonlinear KdV equation. According to the result given in \cite[Proposition 3.1]{MR2463799}, one should expect that for such equations involving space derivatives until the order $3$, the cost of fast controls is bounded by $Ce^{\frac{C}{\sqrt{T}}}$ because of the weights used in the Carleman estimates. We give a precise estimate of $C$ and prove that this power of $T$ is optimal.  

The cost of fast controls for anomalous diffusion equations has been studied notably in \cite{MR2272076} and \cite{MR2248170}, the results are improved in \cite[Section 4.1]{MR2679651}, the latter article gives the optimal power of $1/T$ involved (see also \cite{MR2859866}), but the techniques (spectral inequalities and the Lebeau-Robbiano method) are very different from what we are going to do in this article. In all these articles, the authors were interested in distributed controls in a (small) open subset of the space domain. We consider here only boundary controls, which implies that we will have to make some restrictions on the powers of the Laplace operator we consider.

Our last example concerns the control of fractional Schrödinger equations. As far as we know, the question of the control of such equations was never studied. As before, we are able to derive a precise upper bound on $C_T$.
 
\subsection{Some definitions and notations}
\begin{definition}
Let $\mathcal C$ be a countable set. We say that a sequence of real numbers $(\lambda_n)_{n\in\mathcal C}$ is \emph{regular} if 
$$\gamma((\lambda_n)_{n\in\mathcal C}):=\underset{m\not =n}{\inf}|\lambda_m-\lambda_n|>0.$$
\end{definition}
From now on, we assume that $B$ is a control of the form 
$$Bu=bu,$$
where $b\in \mathcal D(A)'$ ($U$ is here $\mathbb R$ or $\mathbb C$), and we call
$$b_k=<b,e_k>_{(D(A)',D(A))},$$
where $<,>_{(D(A)',D(A))}$ is here the duality product between $\mathcal D(A)'$ and $\mathcal D(A)$ with pivot space $H$.
It is well-known (see \cite{MR704478} and \cite{MR920808}) that if $||(b_k)_{k\in \mathbb N}||_\infty<+\infty$ and if $(\lambda_k)_{k\geqslant 1}$ is regular, then $B$ is an admissible control operator. From now on, we will always assume that $T$ is small enough (for example $T\in (0,1)$).
In the case where $A$ is positive, our main result is the following:
\begin{theorem}
\label{MTh}
\begin{enumerate}
\item Assume that $(\lambda_n)_{n\geqslant 1}$ is a regular increasing sequence of strictly positive numbers verifying moreover that there exist some $\alpha\geqslant 2$ and some $R>0$ such that 
\begin{gather}\label{puiss}\lambda_n=Rn^\alpha+\underset{n\rightarrow \infty}{O}(n^{\alpha-1}),\end{gather}
and assume that  $b_k\simeq 1$ (in the sense that the sequence $(|b_k|)_{k\in\mathbb N}$ is bounded from below and above by strictly positive constants).
Then system \eqref{sch} is null controllable. Moreover, the cost of the control verifies
\begin{gather*}
C_T\lesssim e^{\frac{K}{(RT)^{1/(\alpha-1)}}},\mbox{  for all }K>\frac{3(2)^{1/(\alpha-1)}\pi^{\alpha/(\alpha-1)}}{4((\sin(\pi/\alpha))^{\alpha/(\alpha-1)}}.
\end{gather*}
\item  Assume that $(\lambda_n)_{n\geqslant 1}$ is a regular increasing sequence of strictly positive numbers verifying moreover that there exists some $\alpha>1$ and some constant $R>0$ such that 
\eqref{puiss} holds.
Assume that  $b_k\simeq 1$.
Then system \eqref{heat} is null controllable.
Moreover, the control can be chosen in the space $C^0([0,T],U)$ and the cost of the control (in norm $L^\infty(0,T)$, so this is also true in $L^2(0,T)$) verifies
\begin{gather*}C_T\lesssim e^{\frac{K}{(RT)^{1/(\alpha-1)}}},\mbox{  for all }K> \frac{3(2)^{1/(\alpha-1)}\pi^{\alpha/(\alpha-1)}}{4((2\sin(\pi/(2\alpha)))^{\alpha/(\alpha-1)})}.\end{gather*}

\end{enumerate}
(the implicit constant in the previous inequalities might depend on $\alpha$ but not on $T$).
\end{theorem}
\begin{remark}In the case $\alpha=2$, we obtain exactly the results of \cite{TT}.
\end{remark}

As we will see, we will need for applications in the dispersive case to consider operators $A$ that are not necessarily positive. In the following theorem, we assume that $A$ is self-adjoint with compact resolvent (but not necessarily positive) with a family of eigenvalues $(\lambda_n)_{n\in\mathbb Z^*}$ verifying that $\lambda_n\rightarrow +\infty$ as $n\rightarrow +\infty$ and $\lambda_n\rightarrow -\infty$ as $n\rightarrow -\infty$, and we consider the corresponding dispersive system \eqref{sch} (of course, the corresponding ``parabolic'' system \eqref{heat} cannot be considered).
\begin{theorem}
\label{MTh-disp}Assume that the sequence of increasing eigenvalues $(\lambda_n)_{n\in\mathbb Z^*}$ of $A$  is a regular sequence of non-zero numbers verifying moreover that there exist some $\alpha>1$ and some constant $R>0$ such that 
\begin{equation}
\begin{aligned}\label{doublepuiss}
\lambda_n&=Rn^\alpha+\underset{n\rightarrow \infty}O(n^{\alpha-1}),&\mbox{ }n>0,
\\\lambda_{-n}&=-Rn^\alpha+\underset{n\rightarrow \infty}O(n^{\alpha-1}), \mbox{ }n<0,
\\ \emph{sgn}(\lambda_n)&=\emph{sgn}(n),
\end{aligned}
\end{equation}
and assume that  $b_k\simeq 1$. Then system \eqref{sch} is null controllable. Moreover, the cost of the control verifies
\begin{gather*}
C_T\lesssim e^{\frac{K}{(RT)^{1/(\alpha-1)}}},\mbox{  for all }K>\frac{3(2)^{(\alpha+1)/(\alpha-1)}\pi^{\alpha/(\alpha-1)}}{4((\sin(\pi/\alpha))^{\alpha/(\alpha-1)}}.
\end{gather*}
(the implicit constant in the previous inequalities might depend on $\alpha$ but not on $T$)

\end{theorem}
We are also going to prove that the power of $1/T$ involved in the expression of the cost is optimal in the following sense:
\begin{theorem}\label{lower-bound}
With the same notations and under the same hypothesis as in Theorems~\ref{MTh} and~\ref{MTh-disp}, the power of $1/T$ involved in the exponential is optimal, in the sense that there exists some constant $C>0$ such that in both cases of \eqref{heat} and \eqref{sch} one has  
\begin{gather}\label{low}e^{\frac{C}{T^{1/(\alpha-1)}}}\lesssim C_T.\end{gather}
(the implicit constant in the previous inequality might depend on $\alpha$ but not on $T$)
\end{theorem}
\section{Proofs of Theorems~\ref{MTh}-\ref{lower-bound}}
\subsection{Proof of Theorem~\ref{MTh}}
The following lemma is a refinement of the estimates proved in \cite[Lemma 3.1]{71FR} and is strongly inspired by \cite[Lemma 4.1]{TT}.
\begin{lemma}
\label{gener}Let $(\lambda_n)_{n\geqslant 1}$ be a regular increasing sequence of strictly positive numbers verifying moreover that there exists some $\alpha\geqslant 2$ and some constant $R>0$ such that \eqref{puiss} holds.

Let $\Phi_n$  be defined as follows:
$$\Phi_n(z):=\prod_{k\not =n}(1-\frac{z}{\lambda_k-\lambda_n}).$$
Then 
\begin{enumerate}
\item If $z\in \mathbb C$, \begin{gather}\label{gener-disp}\Phi_n(z)\lesssim e^{\frac{\pi}{\sqrt{R}\sin(\pi/\alpha)}|z|^{\frac{1}{\alpha}}}P(|z|),\end{gather}
where $P$ is a polynomial.

\item If $x\in\mathbb R$,
\begin{gather} \label{gener-parab}\Phi_n(-ix-\lambda_n)\lesssim e^{\frac{\pi}{2\sqrt{R}\sin(\pi/2\alpha)}|x|^{\frac{1}{\alpha}}}\overline P(\lambda_n,|x|),\end{gather}
where $\overline P$ is a polynomial.

(In the previous inequalities, the implicit constant may depend on $\alpha$ but not on $z$, $x$ or $n$)
\end{enumerate}
\end{lemma}
\begin{remark} \label{opt-sup-2}One can see numerically that inequalities \eqref{gener-disp} and \eqref{gener-parab} are optimal for $\alpha\geqslant 2$, but are false for $\alpha\in (1,2)$ (but one could find a less precise estimate).
\end{remark}
\textbf{Proof of Lemma~\ref{gener}.}
Without loss of generality, we can assume that $R=1$ (one can go back to the general case by an easy scaling argument). We have then the existence of some constant $C>0$ such that 
$|\lambda_n-n^\alpha|\leqslant Cn^{\alpha-1}$. From now on we call $\gamma:=\gamma((\lambda_n)_{n\geqslant 1})$.
As in \cite[Page 81]{TT}, one has 
\begin{gather}\label{est-phi-ln}\mbox{ln}|\Phi_n(z)|\leqslant \int_0^{|z|}\int_{\gamma}^\infty \frac{L_n(s)}{(t+s)^2}dsdt,\end{gather}
where $$L_n(s):=\#\{k||\lambda_k-\lambda_n|\leqslant s\}.$$
Let us estimate precisely $L_n(s)$.

One has $$|\lambda_k-\lambda_n|\leqslant s$$ if and only if \begin{gather}\label{cond-sup}\lambda_k-\lambda_n\leqslant s\end{gather} and \begin{gather}\label{cond-inf}\lambda_n-\lambda_k\leqslant s.\end{gather}
\begin{enumerate}
\item Assume that \eqref{cond-sup} holds. Then
$$k^{\alpha-1}(k-C)\leqslant \lambda_n+s.$$
Let $$\overline R(X)=X^{\alpha-1}(X-C).$$
We call $D=\lambda_n+s$.
By studying function $\overline R$, we see that $\overline R(0)=0$, $\overline R(+\infty)=+\infty$ and that $\overline R$ is strictly decreasing on $[0,C(1-1/\alpha)]$ and then strictly increasing on $[C(1-1/\alpha),\infty)$. Hence the equation $\overline R(X)=D$ has a unique solution $\tilde X$ for $n$ sufficiently large and the inequality $\overline R(X)\leqslant D$ is equivalent to $0\leqslant X\leqslant \tilde X$. Moreover, 
$$\overline R(D^{\frac{1}{\alpha}})-D=-CD^{\frac{\alpha-1}{\alpha}}<0$$ and $$\overline R(D^{\frac{1}{\alpha}}+C)-D=(D^{\frac{1}{\alpha}}+C)^{\alpha-1}D^\frac{1}{\alpha}-D=D((1+CD^{-\frac{1}{\alpha}})^{\alpha-1}-1)>0.$$ So $\tilde X\in [D^{1/\alpha},D^{1/\alpha}+C]$ and  $$0\leqslant k\leqslant \tilde X$$ implies \begin{gather}\label{Ln1}k\leqslant (\lambda_n+s)^{\frac{1}{\alpha}}+C.\end{gather}
\item Assume now that \eqref{cond-inf} holds.
\begin{gather}\label{est-t}\lambda_n-s\leqslant k^{\alpha-1}(k+C).\end{gather}
We call $E=\lambda_n-s$. If $\lambda_n-s<0$ then inequality \eqref{est-t} is always true. If $\lambda_n+s\geqslant 0$, we introduce 
$$\tilde R(X)=X^{\alpha-1}(X+C).$$
By studying function $\tilde R$, we see that $\tilde R(0)=0$, $\tilde R(+\infty)=+\infty$ and that $\tilde R$ is strictly increasing on $[0,\infty)$. Hence the equation $\tilde R(X)=E$ has a unique solution $\bar X\in [0,\infty)$ and the inequality $\tilde R(X)\geqslant D$ is equivalent to $ X\geqslant \bar X$. Moreover, $$\tilde R(E^{\frac{1}{\alpha}})-E=CE^{\frac{\alpha-1}{\alpha}}>0$$ and $$\tilde R((E^{\frac{1}{\alpha}}-C)^+)-E=((E^{\frac{1}{\alpha}}-C)^+)^{\alpha-1}E^{\frac{1}{\alpha}}-E=E(((1-CE^{-\frac{1}{\alpha}})^+)^{\alpha-1}-1)<0.$$ So $$\bar X\in [E^{1/\alpha}-C,E^{1/\alpha}]$$ and  $k\geqslant \tilde X$ implies \begin{gather}\label{Ln2}k\geqslant ((\lambda_n-s)^{\frac{1}{\alpha}}-C)^+\geqslant ((\lambda_n-s)^{\frac{1}{\alpha}}-C).\end{gather}
\end{enumerate}
Finally, if we have simultaneously the conditions \eqref{cond-sup} and \eqref{cond-inf}, then combining inequalities \eqref{Ln1} and \eqref{Ln2} necessarily $$k\in [((\lambda_n-s)^{+})^{\frac{1}{\alpha}}-C,(\lambda_n+s)^{\frac{1}{\alpha}}+C]$$ and 
\begin{gather} \label{est-Ln}L_n(s)\leqslant (\lambda_n+s)^{\frac{1}{\alpha}}-((\lambda_n-s)^{+})^{\frac{1}{\alpha}}+2C.\end{gather}
Finally, from \eqref{est-phi-ln} and \eqref{est-Ln},
\begin{gather}\label{Phiest}|\Phi_n(z)|\lesssim (1+|z|/\gamma)^{2C}e^{\int_0^{|z|}\int_{\gamma}^\infty \frac{(\lambda_n+s)^{\frac{1}{\alpha}}-((\lambda_n-s)^{+})^{\frac{1}{\alpha}}}{(t+s)^2}dsdt}.\end{gather}

One has (using the change of variables $v=s/\lambda_n$ for the last inequality)
\begin{equation}\label{PhiUV}
\begin{aligned}\int_0^{|z|}\int_{\gamma}^\infty \frac{(\lambda_n+s)^{\frac{1}{\alpha}}-((\lambda_n-s)^{+})^{\frac{1}{\alpha}}}{(t+s)^2}dsdt&\leqslant |z|\int_{\gamma}^\infty \frac{(\lambda_n+s)^{\frac{1}{\alpha}}-((\lambda_n-s)^{+})^{\frac{1}{\alpha}}}{s(s+|z|)}ds\\ &\leqslant \frac{|z|}{\lambda_n^{1-\frac{1}{\alpha}}}(U(\frac{|z|}{\lambda_n})+V(\frac{|z|}{\lambda_n})),\end{aligned}
\end{equation}
where \begin{gather}\label{U}U(x):=\int_0^1 \frac{(1+v)^{\frac{1}{\alpha}}-(1-v)^{\frac{1}{\alpha}}}{v(v+x)}dv\end{gather}
and \begin{gather}\label{V}V(x):=\int_1^\infty \frac{(v+1)^{\frac{1}{\alpha}}}{v(v+x)}dv.\end{gather}
To prove inequality \eqref{gener-disp}, in view of \eqref{Phiest} and \eqref{PhiUV} it is now enough to prove 
\begin{gather}\label{UetV}x^{1-\frac{1}{\alpha}}(U(x)+V(x))\leqslant \frac{\pi}{\sin(\frac{\pi}{\alpha})}\end{gather}
for all $x\geqslant 0$.

Let us now prove inequality \eqref{UetV}. 
 Let us first study $x^{1-1/\alpha}V(x)$. We remark that 
$$x^{1-1/\alpha}V(x)=x^{1-1/\alpha}\int_1^\infty \frac{(v+1)^{\frac{1}{\alpha}}}{v(v+x)}dv=\int_1^\infty \frac{(v/x+1/x)^{\frac{1}{\alpha}}}{v(v/x+1)}dv.$$
By considering the change of variables $t=x/v$, we obtain

\begin{gather}\label{exprV}x^{1-1/\alpha}V(x)=\int_0^x \frac{(1/t+1/x)^{\frac{1}{\alpha}}}{1+t}dt.\end{gather}
Similarly one has 
\begin{gather}\label{exprU}x^{1-1/\alpha}U(x)=\int_x^\infty \frac{(1/t+1/x)^{\frac{1}{\alpha}}-(1/x-1/t)^{\frac{1}{\alpha}}}{1+t}dt.\end{gather}
Using the dominated convergence Theorem, one proves easily that  
$$x^{1-1/\alpha}V(x)\underset{x\rightarrow \infty}\rightarrow\int_0^\infty \frac{dt}{t^{\frac{1}{\alpha}}(1+t)}$$
and 
$$x^{1-1/\alpha}U(x)\underset{x\rightarrow \infty}\rightarrow 0.$$
Let us call 
$$I(\alpha):=\int_0^\infty \frac{dt}{t^{\frac{1}{\alpha}}(1+t)}.$$
One can compute explicitly this integral.
\begin{lemma}\label{hypergeom}
Let $x>1$. Then 
$$I(x)=\frac{\pi}{\sin(\pi/x)}.$$
\end{lemma}
\textbf{Proof of Lemma~\ref{hypergeom}.}
We remind the following Definition of the Euler Beta function $B$ (see \cite[Page 142, 5.12.3]{MR2723248}):
$$B(x,y):=\int_0^\infty\frac{t^{x-1}}{(1+t)^{x+y}}dt.$$
We then have 
\begin{gather}\label{IBeta}I(x)=B(1-1/x,1/x).\end{gather}
Using the link between the $B$ function and the $\Gamma$ function, we obtain
\begin{gather}\label{eq-gamma}B(1-1/x,1/x)=\frac{\Gamma(1-1/x)\Gamma(1/x)}{\Gamma(1-1/x+1/x)}=\Gamma(1-1/x)\Gamma(1/x).\end{gather}
Using the Euler reflection formula (which can be applied here because $1/x\in(0,1)$), we obtain the desired result.
\cqfd
We will prove that for all $x>0$ one has 
\begin{gather}\label{compI}x^{1-\frac{1}{\alpha}}(U(x)+V(x))\leqslant I(\alpha).\end{gather}
Let us remark that one can compute explicitly $V$ in terms of linear combining of hypergeometric functions: one can use for example Mathematica to check that 
\begin{equation}\label{valV}
\begin{aligned}x^{1-1/\alpha}V(x)&=-\alpha x^{-1/\alpha}\mbox{ }_2F_1(-1/\alpha,-1/\alpha,1-1/\alpha,-1)\\\mbox{ }&+\alpha (1+1/x)^{1/\alpha}\mbox{ }_2F_1(-1/\alpha,-1/\alpha,1-1/\alpha,(x-1)/(x+1)),\end{aligned} 
\end{equation}
where $\mbox{ }_2F_1$ is the ordinary hypergeometric function.
It is then easy to prove that for all $\alpha\geqslant 2$, $x\mapsto x^{1-1/\alpha}V$ is increasing by differentiating \eqref{valV} with respect to $x$.
Let us consider two different cases:
\begin{enumerate}
\item Assume $x<1$. In this case, 
\begin{gather}\label{V1}x^{1-1/\alpha}V(x)\leqslant -\alpha_2F_1(-1/\alpha,-1/\alpha,1-1/\alpha,-1)+\alpha 2^{1/\alpha}.\end{gather}
We remark (by differentiating $x^{1-1/\alpha}U(x)$ with respect to $\alpha$ in expression \eqref{exprU}) that $\alpha\mapsto x^{1-1/\alpha}U(x)$ is increasing, so 
\begin{gather}\label{U1}x^{1-1\alpha}U(x)\leqslant \sqrt{x}\int_0^1 \frac{(1+v)^{\frac{1}{2}}-(1-v)^{\frac{1}{2}}}{v(v+x)}dv\leqslant 1.\end{gather}
(the last inequality in \eqref{U1} can be checked numerically for $x\in[0,1]$)

We also have (the function $\alpha\mapsto\mbox{ }_2F_1(-1/\alpha,-1/\alpha,1-1/\alpha,-1)$ is increasing)
\begin{gather}\label{V2}-\alpha_2F_1(-1/\alpha,-1/\alpha,1-1/\alpha,-1)\leqslant -\alpha_2F_1(-1/2,-1/2,1-1/2,-1)\leqslant-0.52\alpha.\end{gather}

Combining \eqref{V1}, \eqref{U1} and \eqref{V2}, we deduce
$$x^{1-1/\alpha}(U(x)+V(x))\leqslant 1+\alpha 2^{1/\alpha}-0.52\alpha.$$
We just have to prove that 
\begin{gather}\label{okineg}1-0.52\alpha+\alpha 2^{1/\alpha}\leqslant \frac{\pi}{\sin(\pi/\alpha)}.\end{gather}
One verifies numerically that \eqref{okineg} it is true for $\alpha\in[2,3]$,  and one verifies easily by differentiating the expression with respect to $\alpha$ that $\alpha\mapsto 1-0.52\alpha+\alpha 2^{1/\alpha}-\frac{\pi}{\sin(\pi/\alpha)}$ is decreasing at least on $(3,\infty)$.
Inequality \eqref{compI} is proved at least for $x<1$. 
\item Assume $x\geqslant 1$. We have (the equality can be easily obtained thanks to Mathematica for example)
\begin{equation}\label{majU}
\begin{aligned}x^{1-1/\alpha}U(x)&\leqslant x^{-1/\alpha}\int_0^1 \frac{(1+v)^{\frac{1}{\alpha}}-(1-v)^{\frac{1}{\alpha}}}{v}dv\\ \mbox{ }&= x^{-1/\alpha}(H_{1/\alpha}+\mbox{ }_2F_1(-1/\alpha,-1/\alpha,1-1/\alpha,-1)),\end{aligned} 
\end{equation}
where we call $H_{1/\alpha}$ the (generalized) harmonic number of order $1/\alpha$. We have 
\begin{gather}\label{inegH}H_{1/\alpha}\leqslant H_{1/2}\leqslant 0.62.\end{gather}
Using \eqref{valV}, \eqref{majU} and \eqref{inegH}, we deduce
\begin{gather}\label{finUetV}x^{1-1/\alpha}(U(x)+V(x))\leqslant x^{-1/\alpha} A(\alpha)+B(x)\end{gather} with 
$$A(\alpha)=(0.62-\alpha)_2F_1(-1/\alpha,-1/\alpha,1-1/\alpha,-1)$$ and 
$$B(x)=\alpha (1+1/x)^{1/\alpha}\mbox{ }_2F_1(-1/\alpha,-1/\alpha,1-1/\alpha,(x-1)/(x+1)).$$
\end{enumerate}
One has $A(\alpha)<0$, moreover, one easily proves that $B$ is increasing with respect to $x$ and tends to $\alpha_2F_1(-1/\alpha,-1/\alpha,1-1/\alpha,1)=I(\alpha)$.
Hence inequality \eqref{finUetV} implies that inequality \eqref{compI} is also proved for $x\geqslant 1$ and finally 
\eqref{UetV} is proved.
\cqfd

Inequality \eqref{gener-parab} is easier to prove. Doing as in \cite[Page 83]{TT}, we have 
\begin{gather}\label{exprphin}|\Phi_n(-ix-\lambda_n)|^2=\prod_{k\not =n}\frac{|1+ix/\lambda_k|^2}{(1-\lambda_n/\lambda_k)^2}=B_n^2\prod_{k\not =n}|1+x^2/\lambda_k^2|\end{gather}
where 
$$B_n:=\prod_{k\not =n}(1-\lambda_n/\lambda_k)^{-1}.$$
Let us remark that 
\begin{gather}\label{lnOK}\sum_{k\geqslant 1}\ln(1+x^2/\lambda_k^2)= \int_0^{|x|^2/\lambda_1^2} \frac{M(t)}{1+t}dt,\end{gather}
where $$M(t):=\sum_{\lambda_k\leqslant |x|/\sqrt{t}} 1.$$
One easily observe using same computations as before that 

\begin{gather}\label{MOK}M(t)\leqslant |x|^{\frac{1}{\alpha}}t^{-1/(2\alpha)}+C.\end{gather}

We then obtain, using \eqref{lnOK} and \eqref{MOK},
$$\sum_{k\geqslant 1}\ln(1+x^2/\lambda_k^2) \leqslant C\ln(1+|x|^2/\lambda_1^2)+|x|^{\frac{1}{\alpha}} \int_0^{\infty} \frac{1}{t^{1/(2\alpha)}(1+t)}dt\leqslant C\ln(1+|x|^2/\lambda_1^2)+|x|^{\frac{1}{\alpha}}I(2\alpha).$$
We deduce by Lemma~\ref{hypergeom} and \eqref{exprphin} that
$$\Phi_n(-ix-\lambda_n)\lesssim B_n (1+|x|^2/\lambda_1^2)^{C/2}e^{\pi |x|^{\frac{1}{\alpha}}/(2\sin(\pi/(2\alpha)))}$$
and it can be proved that $B_n$ is at most polynomial in $\lambda_n$ (the computations are the same as in \cite[Pages 83-84]{TT}) as wished. This proves inequality \eqref{gener-parab}.
\cqfd

Now, we study the multiplier, which is very similar to the one studied in \cite{TT}. Let $\nu>0$ and $\beta>0$ be linked by the following relation: 
\begin{gather}\label{betanu}\beta \nu^{\alpha-1}=((\pi+\delta)/(\sin(\pi/\alpha)))^\alpha,\end{gather}
where $\delta>0$ is a small parameter.

We call
$$\sigma_\nu(t):=\exp(-\frac{\nu}{(1-t^2)})$$ prolonged by $0$ outside $(-1;1)$. $\sigma_\nu$ is analytic on $B(0,1)$. We call
$$H_\beta(z):=C_\nu\int_{-1}^1\sigma_\nu(t)e^{-i\beta t z}dt,$$
where 
$$C_\nu:=1/||\sigma_\nu||_1.$$
Thanks to \cite[Lemma 4.3]{TT}, we have 
\begin{align}\label{H0}H_\beta(0)&=1,
\\
\label{minmult}H_\beta(ix)&\gtrsim \frac{e^{\beta|x|/(2\sqrt{\nu+1})}}{\sqrt{\nu+1}},\\
\label{maj-Cv} \frac{1}{2}e^\nu&\leqslant C_\nu\leqslant \frac{3}{2}\sqrt{\nu+1}e^\nu,\\
\label{paley-wiener}|H_\beta(z)|&\leqslant e^{\beta|Im(z)|}.
\end{align}
The main estimate is the following:
\begin{lemma}
\label{est-mul}For $x\in\mathbb R$, we have
$$H_{\beta}(x)\lesssim \sqrt{\nu+1} e^{3\nu/4-((\pi+\delta/2)|x|^{\frac{1}{\alpha}})/(\sin(\pi/\alpha))}.$$
(The implicit constant may depend on $\alpha$)
\end{lemma}
\begin{remark}\label{inf-2}Lemma~\ref{est-mul} is false for $\alpha\in (1,2)$. This explains why we were not able to extend Theorem~\ref{MTh} to the case where $\alpha\in (1,2)$. However, we know that systems like \eqref{heat} and \eqref{sch} are null controllable as soon as $\alpha>1$, so one can conjecture that there is a way to extend the previous estimates for $\alpha\in(1,2)$.
\end{remark}
\textbf{Proof of Lemma~\ref{est-mul}.}
First of all, consider some $t\in[0,1)$ and $\theta\in(-\pi,\pi)$. We call $\rho:=1-t$ and $z:=t+\rho e^{i\theta}$. One has (see \cite[Page 85]{TT})
$$Re\frac{1}{1-z^2}\geqslant \frac{1}{4\rho}+\frac{1}{4}\geqslant \frac{1}{4\rho^{1/(\alpha-1)}}+\frac{1}{4},$$
because $\rho\leqslant 1$ and $\alpha\geqslant 2$.
So, doing as in \cite{TT}, we obtain by applying the Cauchy formula for holomorphic functions
$$|\sigma_\nu^{(j)}(t)|\leqslant j! e^{-\frac{\nu}{4}}\underset{\rho>0}{\sup} \frac{e^{{-\frac{1}{4\rho^{1/(\alpha-1)}}}}}{\rho^j}.$$
Computing the supremum on $\rho\in \mathbb R^{+*}$, we obtain 
\begin{gather}\label{infsigma}|\sigma_\nu^{(j)}(t)|\leqslant j! e^{-\frac{\nu}{4}}
e^{-(\alpha-1)j}(4^{1-\alpha}((\alpha-1)j)^{1-\alpha})^{-j},\mbox{ }t\in [0,1).
\end{gather}
Using the fact that $\sigma_\nu$ is even, inequality \eqref{infsigma} is true for every $t\in(-1,1)$.
Using inequality $j!>j^je^{-j}$ in \eqref{infsigma}, we obtain 
\begin{gather}\label{maj-sj}|\sigma_\nu^{(j)}(t)|\leqslant (j!)^{\alpha} e^{-\frac{\nu}{4}}
(4^{1-\alpha}(\alpha-1))^{-j}.\end{gather}
Since all derivatives of $\sigma_\nu$ vanish at $t=-1$ and $t=1$, we have 
\begin{gather}\label{cool-est}|H_\beta(x)|\leqslant \frac{2C_\nu ||\sigma_\nu^{(j)}||_\infty}{(\beta x)^j},\end{gather}
for all $x>0$ and $j\in\mathbb N$. Combining \eqref{maj-sj}, \eqref{cool-est} and \eqref{maj-Cv}, we deduce that 
\begin{gather}\label{Hb-j}|H_\beta(x)|\lesssim \sqrt{\nu+1}(j!)^\alpha e^{\frac{3\nu}{4}}\frac{4^{(\alpha-1)j}}{((\alpha-1)\beta x)^j},\mbox{ }j\in\mathbb N.\end{gather}
We set \begin{gather}\label{bonj} j:=\lfloor (1/a)((\alpha-1)\beta x)^{1/\gamma}\rfloor\end{gather} with some constants $a$ and $\gamma$ which will be chosen correctly soon.
Then we have \begin{gather}\label{inegj}(\alpha-1)\beta x\geqslant (aj)^{\gamma}.\end{gather} Using \eqref{inegj} and  \eqref{Hb-j} we obtain 
\begin{gather}\label{Hf}|H_\beta(x)|\lesssim \sqrt{\nu+1}(j!)^\alpha e^{\frac{3\nu}{4}}\frac{4^{(\alpha-1)j}}{(aj)^{\gamma j}}.\end{gather}
We choose $\gamma=\alpha$ and  $a=4^{1-1/\alpha}$. Combining \eqref{Hf}, \eqref{bonj}, \eqref{betanu} and inequality 
\begin{gather*}\label{estej}(j!)^\alpha \leqslant j^{\alpha/2} j^{\alpha j} e^{-\alpha j},\end{gather*}
we deduce 
$$|H_\beta(x)|\lesssim  \sqrt{\nu+1}e^{\frac{3\nu}{4}}e^{-\alpha j}j^{\alpha/2}\leqslant \sqrt{\nu+1}e^{\frac{3\nu}{4}}e^{-(\pi +\delta/2)/(sin(\pi/\alpha))|x|^{\frac{1}{\alpha}}}.$$
\cqfd
\textbf{Proof of Theorem~\ref{MTh}.}

The proof follows the proof of \cite[Theorem 3.1 and 3.4]{TT}. We still assume without loss of generality that $R=1$. Let us first consider the dispersive case (Equation \eqref{sch}).
We call 
\begin{gather}\label{def-gn}g_n(z):=\Phi_n(-z-\lambda_n)H_\beta(z+\lambda_n).\end{gather}
We want to apply at the end the Paley-Wiener Theorem (see estimate \eqref{paley-wiener}) in an optimal way, so we want $\beta$ to be close to $T/2$. Assume that $\beta<T/2$ and close to $T/2$, for example \begin{gather}\label{beta-d}\beta=\frac{T(1-\delta)}{2}.\end{gather}
One has $g_n(-\lambda_k)=\delta_{kn}$. Moreover, thanks to \eqref{def-gn}, \eqref{gener-disp}, Lemma~\ref{est-mul}, \eqref{betanu} and \eqref{beta-d}
\begin{align*}|g_n(x)|&\lesssim e^{\frac{3\nu}{4}+\pi/\sin(\pi/\alpha)|x+\lambda_n|^{\frac{1}{\alpha}}-(\pi+\delta/2)/\sin(\pi/\alpha)|x+\lambda_n|^{\frac{1}{\alpha}}}P(|x+\lambda_n|) \\ 
&\lesssim e^{\frac{3\nu}{4}-\delta/(2\sin(\pi/\alpha))|x+\lambda_n|^{\frac{1}{\alpha}}}P(|x+\lambda_n|)\\&\lesssim \frac{e^{\frac{3}{4}(\pi+\delta)^{\alpha/(\alpha-1)}/((\sin(\pi/\alpha))^{\alpha/(\alpha-1)}\beta^{1/(\alpha-1)})}}{1+(x+\lambda_n)^2}
\\ &\lesssim \frac{e^{\frac{3}{4}2^{1/(\alpha-1)}(\pi+\delta)^{\alpha/(\alpha-1)}/((\sin(\pi/\alpha))^{\alpha/(\alpha-1)}(T(1-\delta))^{1/(\alpha-1)})}}{1+(x+\lambda_n)^2}.\end{align*}
Considering $\delta$ as close as $0$ as needed, we deduce that  
\begin{gather}\label{estgn}|g_n(x)|\lesssim \frac{e^{\frac{K}{T^{1/(\alpha-1)}}}}{1+(x+\lambda_n)^2}\end{gather}
for all $$K>\frac{3}{4}2^{1/(\alpha-1)}\pi^{\alpha/(\alpha-1)}/(\sin(\pi/\alpha))^{\alpha/(\alpha-1)}.$$
This notably proves that $g_n\in L^2(\mathbb R)$. Moreover, using \eqref{gener-disp}, \eqref{def-gn}, \eqref{beta-d} and \eqref{paley-wiener}, we obtain
$$|g_n(z)|\lesssim e^{T|z|/2}.$$
Hence, using the Paley-Wiener Theorem, $g_n$ is the Fourier transform of a function $f_n\in L^2(\mathbb R)$ with compact support $[-T/2,T/2]$. Moreover, by construction $\{f_n\}$ is biorthogonal to the family $\{e^{i\lambda_nt}\}$.
Then, one can create the control thanks to the family $\{f_n\}$. Let us consider $y^0=\sum a_k e_k$ the initial condition, we call 
\begin{gather}\label{defud}u(t):=-\sum_{k\in\mathbb N} (a_k/b_k)e^{-iT\lambda_k/2}f_k(t-T/2).\end{gather}
This expression is meaningful since $b_k\simeq 1$, moreover
the corresponding solution $y$ of \eqref{sch} verifies $y(T,\cdot)\equiv 0$.
 Using the Minkovski inequality, Parseval equality, \eqref{defud}, $b_k\simeq 1$ and \eqref{estgn}, we obtain
\begin{align*}||u(t)||_{L^2(0,T)}&\lesssim e^{\frac{K}{T^{1/(\alpha-1)}}}(\sum |a_k|^2(\int_{\mathbb R} \frac{dx}{(1+(x+\lambda_n)^2)^2}))^{1/2}\\&\lesssim e^{\frac{K}{T^{1/(\alpha-1)}}}(\pi/2\sum |a_k|^2)^{1/2}\\&\lesssim e^{\frac{K}{T^{1/(\alpha-1)}}}||y^0||_{H}.\end{align*}

We now consider the parabolic case (Equation \eqref{heat}).
We call 
\begin{gather}\label{defhn}h_n(z):=\frac{\Phi_n(-iz-\lambda_n)H_\beta(z \sin(\pi/\alpha)^\alpha/(2\sin(\pi/(2\alpha))^\alpha))}{H_\beta(i\lambda_n\sin(\pi/\alpha)^\alpha/(2\sin(\pi/(2\alpha))^\alpha))}.\end{gather}
 Assume that $$\beta<\frac{T(2\sin(\pi/2\alpha))^\alpha}{2\sin(\pi/(\alpha))^\alpha}$$ and close to $$\frac{T (2\sin(\pi/2\alpha))^\alpha}{2(\sin(\pi/(\alpha)))^\alpha},$$ for example \begin{gather}\label{betah}\beta=\frac{(1-\delta)T (2\sin(\pi/2\alpha))^\alpha}{2\sin(\pi/(\alpha))^\alpha}.\end{gather}
One has $h_n(i\lambda_k)=\delta_{kn}$. Moreover, thanks to \eqref{defhn}, \eqref{gener-parab}, \eqref{minmult}, Lemma~\ref{est-mul}, \eqref{betanu} and \eqref{betah}, one has 
\begin{align*}|h_n(x)|\lesssim&(\nu +1) e^{\frac{3}{4}\nu+\pi/(2\sin(\pi/2\alpha))|x|^{\frac{1}{\alpha}}-((\pi+\delta/2)/(2\sin(\pi/2\alpha)))|x|^{\frac{1}{\alpha}}-\frac{\beta|\lambda_n|}{2\sqrt{\nu+1}}}\overline P(|x|,|\lambda_n|) \\ &
\lesssim (\nu +1) e^{\frac{3}{4}\nu-\delta/(2\sin(\pi/2\alpha))|x|^{\frac{1}{\alpha}}-\frac{\beta\lambda_n}{2\sqrt{\nu+1}}}P(|x|,\lambda_n|)\\&\lesssim (\nu +1)\frac{e^{\frac{3}{4}2^{1/(\alpha-1)}(\pi+\delta)^{\alpha/(\alpha-1)}/((2\sin(\pi/\alpha))^{\alpha/(\alpha-1)}\beta^{1/(\alpha-1)})}}{(1+(x+\lambda_n)^2)}
\\&\lesssim (\nu +1)\frac{e^{\frac{3}{4}2^{1/(\alpha-1)}(\pi+\delta)^{\alpha/(\alpha-1)}/((2\sin(\pi/(2\alpha)))^{\alpha/(\alpha-1)}(T(1-\delta))^{1/(\alpha-1)})}}{(1+(x+\lambda_n)^2)}.\end{align*}
Considering $\delta$ as close as $0$ as needed, we deduce that  
\begin{gather}\label{esthn}|h_n(x)|\lesssim \frac{e^{\frac{K}{T^{1/(\alpha-1)}}}}{(1+(x+\lambda_n)^2)},\end{gather}
for all $$K>\frac{3}{4}2^{1/(\alpha-1)}\pi^{\alpha/(\alpha-1)}/((2\sin(\pi/(2\alpha)))^{\alpha/(\alpha-1)}).$$
This notably implies that $h_n(x)\in L^1(\mathbb R)\cap L^2(\mathbb R)$ and 
\begin{gather}\label{inttf} ||h_n||_{L^1(\mathbb R)}\lesssim e^{\frac{K}{T^{1/(\alpha-1)}}}.\end{gather} Moreover, using  \eqref{gener-parab}, \eqref{defhn}, \eqref{paley-wiener} and \eqref{betah}
$$|h_n(z)|\lesssim e^{T|z|/2},$$
so using the Paley-Wiener Theorem, $h_n$ is the Fourier transform of a function $w_n\in L^2(\mathbb R)$ with compact support $[-T/2,T/2]$. Moreover, by construction $\{w_n\}$ is biorthogonal to the family $\{e^{-\lambda_nt}\}$.
Then, one can create the control thanks to the family $\{h_n\}$. Let us consider $y^0=\sum a_k e_k$ the initial condition, we call 
\begin{gather}\label{defuh}u(t):=-\sum (a_k/b_k)e^{-T\lambda_k/2}w_k(t-T/2),\end{gather}
This expression is meaningful since $b_k\simeq 1$, moreover
the corresponding solution $y$ of \eqref{heat} verifies $y(T,\cdot)\equiv 0$. One easily verifies that $u\in C^0([0,T],\mathbb R)$. Using \eqref{defuh}, $|b_k|\simeq 1$ and inequality \eqref{inttf}, we obtain
$$||u(t)||_{L^\infty(0,T)}\lesssim e^{\frac{K}{T^{1/(\alpha-1)}}}\sum |a_k|e^{-T\lambda_k/2}.$$
Using the Cauchy-Schwarz inequality, one deduces
$$||u(t)||_{L^\infty(0,T)}\lesssim e^{\frac{K}{T^{1/(\alpha-1)}}}||y^0||_H.$$
\cqfd
\subsection{Proof of Theorem~\ref{MTh-disp}}
We will not give the details of the proof of Theorem~\ref{MTh-disp} because it is exactly the same as the one of Theorem~\ref{MTh}. We just explain in details the modifications appearing in Lemma~\ref{gener}.

\begin{lemma}
\label{gener-KdV}Let $(\lambda_n)_{n\in\mathbb Z}$ be a regular increasing sequence of non-zeros numbers verifying moreover that there exists some $\alpha\geqslant 2$ and some constant $R>0$ such that 
\eqref{doublepuiss}
 holds.
Let $\Phi_n$  be defined as follows:
$$\Phi_n(z):=\prod_{k\not =n}(1-\frac{z}{\lambda_k-\lambda_n}),$$
then 
\begin{gather}\label{gener-disp2}\Phi_n(z)\lesssim e^{\frac{2\pi}{\sqrt{R}\sin(\frac{\pi}{\alpha})}|z|^{\frac{1}{\alpha}}}P(|z|),
\end{gather}
where $P$ is a polynomial in $|z|$. 
(In the previous inequality, the implicit constant may depend on $\alpha$ but not on $z$ or $n$.)
\end{lemma}
\textbf{Proof of Lemma~\ref{gener-KdV}.}
We use the same notations as in the proof of Lemma~\ref{gener} and assume without loss of generality that $R=1$. Let us give a new upper bound for $L_n(s)$.
 
Let $s\geqslant 0$ and let us estimate $\#\{k||\lambda_k-\lambda_n|\leqslant s\}$. If $k$ and $n$ have the same sign, we have necessarily (see the proof of Lemma~\ref{gener} and \eqref{doublepuiss})
\begin{gather}\label{estk1}\#\{k||\lambda_k-\lambda_n|\leqslant s, \emph{sgn}(k)=\emph{sgn}(n)\}\leqslant (|\lambda_n|+s)^{\frac{1}{\alpha}}-((|\lambda_n|-s)^{+})^{\frac{1}{\alpha}}+2C.\end{gather}
If $k$ and $n$ have different sign, one can assume without loss of generality that $k>0$, so that one has $\lambda_k>0$ and $\lambda_n<0$ (see \eqref{doublepuiss}). If 
$|\lambda_k-\lambda_n|\leqslant s$, then necessarily   
$\lambda_k\leqslant (s-|\lambda_n|)^+$, i.e. 
$$k^\alpha-Ck^{\alpha-1}\leqslant D,$$
with $D=s-|\lambda_n|$. Using the same reasoning as in the proof of Lemma~\ref{gener}, this implies that
\begin{gather}\label{estk2}k\leqslant ((|\lambda_n|-s)^{+})^{\frac{1}{\alpha}}+C.\end{gather}
Finally, combining \eqref{estk1} and \eqref{estk2}, we obtain
\begin{gather}\label{newln}L_n(s)\leqslant (|\lambda_n|+s)^{\frac{1}{\alpha}}-((|\lambda_n|-s)^{+})^{\frac{1}{\alpha}}+((s-|\lambda_n|)^{+})^{\frac{1}{\alpha}}+3C.
\end{gather}
We then have using \eqref{est-phi-ln} and \eqref{newln}
$$|\Phi_n(z)|\lesssim (1+|z|/\delta)^{3C}e^{\int_0^{|z|}\int_{\gamma((|\lambda_n|)_{n\geqslant 1})}^\infty \frac{(|\lambda_n|+s)^{\frac{1}{\alpha}}-((|\lambda_n|-s)^{+})^{\frac{1}{\alpha}}+((s-|\lambda_n|)^{+})^{\frac{1}{\alpha}}}{(t+s)^2}dsdt}.$$
One has 
$$\int_0^{|z|}\int_{\gamma}^\infty \frac{(|\lambda_n|+s)^{\frac{1}{\alpha}}-((|\lambda_n|-s)^{+})^{\frac{1}{\alpha}}+((s-|\lambda_n|)^{+})^{\frac{1}{\alpha}}}{(t+s)^2}dsdt\leqslant \frac{|z|}{|\lambda_n|^{1-\frac{1}{\alpha}}}(U(\frac{|z|}{|\lambda_n|})+V(\frac{|z|}{|\lambda_n|})+W(\frac{|z|}{|\lambda_n|})),$$
where $U$ and $V$ have already been defined in \eqref{U} and \eqref{V}, and where
$$W(x):=\int_1^\infty\frac{(u-1)^{1/\alpha}}{u(u+x)}du.$$
Since we already proved by Lemma~\ref{hypergeom} and \eqref{compI} that
$$x^{1-\frac{1}{\alpha}}(U(x)+V(x))\leqslant \frac{\pi}{\sin(\pi/\alpha)},$$
Lemma~\ref{gener-KdV} will be proved as soon as 
\begin{gather}\label{inW}x^{1-\frac{1}{\alpha}}W(x)\leqslant \frac{\pi}{\sin(\pi/\alpha)}.\end{gather}
Using the change of variable $u=1/s$, we obtain 
\begin{gather}\label{est-W1}W(x)=\int_0^1\frac{(1-s)^{1/\alpha}}{s^{1/\alpha}(1+sx)}ds=\int_0^1 \frac{(1-s)^{1/\alpha-1}}{s^{1/\alpha}}\frac{1-s}{1+sx}ds.\end{gather}
One has the equality
\begin{gather}\label{eqfrac}\frac{1-s}{1+sx}=\frac{1+x}{x(1+sx)}-\frac{1}{x}.\end{gather}
Replacing \eqref{eqfrac} in \eqref{est-W1}, we deduce
\begin{gather}\label{est-W2}W(x)=\frac{1+x}{x}\int_0^1\frac{(1-s)^{1/\alpha-1}}{s^{1/\alpha}(1+sx)}ds-\frac{1}{x}\int_0^1\frac{(1-s)^{1/\alpha-1}}{s^{1/\alpha}}ds.\end{gather}
The usual definition of $B$ (see \cite[Page 142, 5.12.1]{MR2723248}) gives 
\begin{gather}\label{defbeta2}\int_0^1\frac{(1-s)^{1/\alpha}}{s^{1/\alpha}}ds=B(1-1/\alpha,1/\alpha).\end{gather}
Using Lemma~\ref{hypergeom},\eqref{compI}, \eqref{est-W2}, \eqref{defbeta2} and the symmetry of the $B$ function, we deduce 
\begin{gather}\label{calc-W}W(x)=\frac{1+x}{x}\int_0^1\frac{(1-s)^{1/\alpha-1}}{s^{1/\alpha}(1+sx)}ds-\frac{\pi}{x\sin(\pi/\alpha)}.\end{gather}
Using the change of variables $u=1/s$, we have 
\begin{gather}\label{chv}\int_0^1\frac{(1-s)^{1/\alpha-1}}{s^{1/\alpha}(1+sx)}dt=\int_1^\infty\frac{(u-1)^{1/\alpha-1}}{(u+x)}du.\end{gather}
Using the change of variables $s=u/(1+x)$, Lemma~\ref{hypergeom},\eqref{compI}, and the symetry of the Beta function, we obtain
\begin{gather}\label{chv2}\int_1^\infty\frac{(u-1)^{1/\alpha-1}}{u+x}du=\int_0^\infty\frac{u^{1/\alpha-1}}{u+1+x}du=(1+x)^{1/\alpha-1}B(1/\alpha,1-1/\alpha)=(1+x)^{1/\alpha-1}\frac{\pi}{\sin(\pi/\alpha)}.\end{gather}
Going back to \eqref{calc-W} and using \eqref{chv} and \eqref{chv2}, we deduce 
$$W(x)=\frac{\pi((x+1)^{\frac{1}{\alpha}}-1)}{x\sin(\pi/\alpha)}.$$
Then, using also the inequality  (true for $x\geqslant 0$ and $\alpha\geqslant 1$)
$$(x+1)^{\frac{1}{\alpha}}-1\leqslant x^{\frac{1}{\alpha}},$$
we obtain \eqref{inW}.
\cqfd
\subsection{Proof of Theorem~\ref{lower-bound}}
 We follow the strategy given in \cite{MR805273}. Without loss of generality we can assume that $R=1$ in \eqref{puiss} and \eqref{doublepuiss}. Looking carefully at this article, one observes that one could adapt the reasoning to equations \eqref{heat} and \eqref{sch}. We treat the case of real or pure imaginary eigenvalues (of $A$ or $iA$, see equations \eqref{heat} and \eqref{sch}) $\lambda_n$ or $i\lambda_n$ with $\lambda_n$ verifying \eqref{puiss} (one could easily adapt the reasoning to obtain the same results in the dispersive case with $\lambda_n$ verifying \eqref{doublepuiss}). We introduce $(\mu_n):=(\lambda_n)$ in the parabolic case and $(\mu_n):=(-i\lambda_n)$ in the dispersive case. We call $$E(T):=\overline{\emph{span}(\{e^{-\mu_nt}|n\in\mathbb N\})}^{L^2(0,T)},$$ $$E_m(T):=\overline{\emph{span}(\{e^{-\mu_nt}|n\not =m\})}^{L^2(0,T)}.$$

We remark that using the results of  \cite{MR0014502} for the parabolic case and \cite{MR0015553} for the dispersive case, if the sequence $(\lambda_n)_{n\in\mathbb N^*}$ verifies \eqref{puiss}, then $E(T)$ in a proper subspace of $L^2(0,T)$ and $e^{-\mu_mt}\not \in E_m(T)$. Moreover, if we call $d_m(T)$ the distance between $e^{-\mu_mt}$ and $E_m(T)$ and $r_m$ the orthogonal projection of $e^{-\mu_mt}$ over $E_m(T)$, then the family $\{\psi_m\}$ defined by 
$$\psi_m(t):=\frac{e^{-\mu_mt}-r_n(t)}{d_m(T)^2}$$ is biorthogonal to the family of exponentials $\{e^{-\mu_mt}\}$ (see \cite{71FR} or \cite{MR805273}, this can be easily generalized in the case of purely imaginary eigenvalues). One also has 
\begin{gather}\label{dmpsim}||\psi_m||_{L^2(0,T)}=\frac{1}{d_m(T)}.\end{gather}
If $y^0:=\sum a_ke_k$, then the control $u$ is given by 
\begin{gather}\label{bestc}u(t):=-\sum a_k/b_k\psi_k(T-t)\end{gather} and one can easily prove that this control is optimal in $L^2(0,T)$. We are now going to give an upper bound on $d_m(T)$, which would provide a lower bound on $C_T$. In all what follows, $C(m)$ denotes some constant depending only on the integer $m$ (and possibly on $\alpha$) that may change from one line to another.
\begin{lemma}\label{gui-est}For every $m\in\mathbb N$, there exists some numerical constant $a(m)$ and some constant $C(m)$ such that 
\begin{gather}\label{estdm}d_m(T)\leqslant C(m)T^{1/2}(j!)^{\alpha-1}(a(m)T)^j\end{gather}
holds for $j\in\mathbb N$ and $T>0$.
\end{lemma}

\begin{remark}As before, we are not able to extend this Lemma to the case where $\alpha\in (1,2)$ (precisely because of estimate \eqref{ineqxaok} which is false in this case), and hence we are not able to extend Theorem~\ref{lower-bound} to this case.
\end{remark}
\textbf{Proof of Lemma~\ref{bok}.}
Following \cite{MR805273}, we only treat the case $j\geqslant m$ (inequality \eqref{estdm} has only an interest for large $j$ because if we prove it for $j\geqslant m$ then it is automatically true for $j<m$ by increasing the constant $C(m)$ in front of the right-hand side if necessary).
One can prove (by considering a finite number of modes, see $(4.9)$ in \cite{MR805273}) that for all $j\geqslant 1$ one has 
\begin{gather}\label{estdm1}d_m(T)\leqslant \frac{T^{j+\frac{1}{2}}}{j!\sqrt{2j+1}}\prod_{r=1}^{m-1}|\lambda_r-\lambda_m|\prod_{r=m}^j|\lambda_{r+1}-\lambda_m|.\end{gather}
For $k,l\in\mathbb N$ one has 
\begin{gather}\label{lkl}|\lambda_k-\lambda_l|\leqslant |k^\alpha-l^\alpha|+C(k^{\alpha-1}+l^{\alpha-1}).\end{gather}

We deduce from \eqref{lkl} that
\begin{equation}\label{lambdamn1}
\begin{aligned}\prod_{r=1}^{m-1}|\lambda_r-\lambda_m|\prod_{r=m}^j|\lambda_{r+1}-\lambda_m|& \mbox{ }\\ \leqslant \prod_{r=1}^{m-1}|r^\alpha-m^\alpha|(1+C\frac{(r^{\alpha-1}+m^{\alpha-1})}{r^\alpha-m^\alpha})\prod_{r=m}^{j}|(r+1)^\alpha-m^\alpha|(1+C\frac{((r+1)^{\alpha-1}+m^{\alpha-1})}{(r+1)^\alpha-m^\alpha}).&\mbox{ }\end{aligned}
\end{equation}
For every $m\in\mathbb N$, there exists some $C(m)>0$ such that for every $j\not=m$,
\begin{gather}\label{est-log}
1+C\frac{(j^{\alpha-1}+m^{\alpha-1})}{j^\alpha-m^\alpha}\leqslant1+\frac{C(m)}{j}.
\end{gather}

Using \eqref{estdm1}, \eqref{lkl},\eqref{lambdamn1} and \eqref{est-log}, we deduce that 
\begin{gather}\label{estdm2}d_m(T)\leqslant C(m)\frac{T^{j+\frac{1}{2}}}{j!\sqrt{2j+1}}\prod_{r=m}^j(1+C(m)/r)|(r+1)^\alpha-m^\alpha|.\end{gather}
One has $\sum_{r=m}^j \ln(1+C(m)/r)\sim C(m)\ln(j)$ as $j\rightarrow \infty$ so 
\begin{gather}\label{prodneg}\prod_{r=m}^j(1+C(m)/r) \lesssim j^{C(m)}.\end{gather} 

Let us study the quantities of the form $k^{\alpha}-l^{\alpha}$ with $k\geqslant l$.
\begin{gather}\label{raminf}|k^{\alpha}-l^{\alpha}|=k^{\alpha}|1-\frac{l^{\alpha}}{k^{\alpha}}|.\end{gather}
One easily verifies that the following inequality holds for $\alpha\geqslant 2$ and $x\in[0,1]$:

\begin{gather}\label{ineqxaok}1-x^\alpha\leqslant (1+x)^{\alpha-1}(1-x).\end{gather}
We deduce from \eqref{raminf} and \eqref{ineqxaok} (the constant $C(m)$ may change from one line to another)
\begin{align*}\mbox{ }&\prod_{r=m}^j|(r+1)^\alpha-m^\alpha|
\\\mbox{ }&\leqslant \prod_{r=m}^j(r+1-m)(r+1+m)^{\alpha-1}
\\\mbox{ }&\leqslant C(m)(j+1-m)!((j+1+m)!)^{\alpha-1}.
\end{align*}

Using the computations above, inequality (true for $j\geqslant m$)
$$(j+1+m)!\leqslant C(m)j^{C(m)}j!,$$
\eqref{estdm2} and \eqref{prodneg}, we deduce that
$$d_m(T)\leqslant C(m)j^{C(m)}(j!)^{\alpha-1} T^{j+\frac{1}{2}},$$
so that \eqref{estdm} holds if we choose $a(m)>0$ large enough such that $j^{C(m)}\leqslant a(m)^j$.
\cqfd

We deduce, using \eqref{dmpsim} and \eqref{estdm}, that for all $j\in\mathbb N$ one has
\begin{gather}\label{finpsim}||\psi_m||_{L^2(0,T)}\geqslant \frac{1}{C(m)T^{1/2}(j!)^{\alpha-1}(a(m)T)^j}.\end{gather}
Using equality \eqref{bestc} with initial condition eigenvector $e_m$ and the following inequality true for $j$ large enough
$$j!\leqslant j^je^{-j/2},$$
and choosing (with $T$ small enough) $$j:=\lceil (1/(a(m)T))^{\frac{1}{\alpha-1}}\rceil,$$ one easily proves using inequality \eqref{finpsim} that \eqref{low} holds.
\cqfd
\section{Applications}
\subsection{Linear KdV equations controlled on the boundary: the case of periodic boundary conditions with a boundary control on the derivative of the state}
In this section, we consider the following controlled linearized KdV equation posed on $(0,T)\times(0,L)$ (this is the first example studied in \cite{97Ro}).
Let us first introduce the following family of periodic Sobolev spaces (endowed with the usual Sobolev norm)
$$H_p^k:=\{y\in H^k(0,L)|u^{(j)}(0)=u^{(j)}(L),j=0\ldots k-1\}.$$
We consider the following equation:
\begin{equation}\label{KdVper}
\begin{aligned}
y_t+y_{xxx}&=0&\mbox{ in }(0,T)\times(0,L),
\\y(t,0)&=y(t,L)&\mbox{ in }(0,T),
\\y_x(t,0)&=y_x(t,L)+u(t)&\mbox{ in }(0,L),
\\y_{xx}(t,0)&=y_{xx}(t,L)&\mbox{ in }(0,L),
\end{aligned}
\end{equation}
with initial condition $y^0\in H:=(H^1_p)'$ and control $u\in L^2(0,T)$.  
This system was first studied in \cite{MR1214759} where the authors proved a result of exact controllability under the technical condition  that the integral in space of the initial state had to be equal to the one of the final state. This result was improved later in \cite{97Ro}.
We know (see \cite[Remark 2.3]{97Ro}) that in this case there exists a unique solution $y\in C^0([0,T],(H_p^1)')$ to \eqref{KdVper}. Moreover, it is explained in \cite[Remark 2.3]{97Ro} that this equation is exactly controllable (and then null controllable) at all time $T>0$ for all length $L>0$ (in fact the case which is treated in \cite{97Ro} is $L=2\pi$ but it can be easily extended to all $L$). We call $\mathcal A$ the operator $\partial^3_{xxx}$ with domain 
$\mathcal D(A):=H^2_p(0,L).$ This operator is skew-adjoint, the  eigenvalues are $i\lambda_k:=8i\pi^3k^3/L^3$ for $k\in\mathbb Z$ , the corresponding eigenfunction is (normed in $(H^1_p)'$) $$e_k:x\mapsto \frac{(1+4\pi^2k^2/L^2)^{1/2}e^{\frac{i2\pi kx}{L}}}{\sqrt{L}}.$$ If $y^0\in (H^1_p)'$ is written under the form $y^0(x)=\sum_{k\in\mathbb Z} a_ke_k(x)$, then 
the solution $y$ of \eqref{KdVper} can be written under the form $$y(t,x)=\sum_{k\in\mathbb Z}a_ke^{i\lambda_k t}e_k(x).$$
One easily proves (using integrations by parts, see for example \cite[Section 2.7, page 101]{MR2302744}) that for every $\varphi\in \mathcal D(A)$,
$$b(\varphi)=-(\Delta^{-1}\varphi)'(0),$$
so that
$$b=\delta_L'\circ \Delta^{-1},$$
where $\Delta^{-1}:=-(-\Delta^{-1})$ is the inverse of the Dirichlet-Laplace operator.
We have $$|b_k|=|e_k'(L)|/k^2\simeq 1.$$
One can apply directly Theorem~\ref{MTh-disp} and Theorem~\ref{lower-bound} with $k=3$ and $R=\frac{8\pi^3}{L^3}$ to obtain:
\begin{theorem}\label{MTh-KdV2}
Equation \eqref{KdVper} is null controllable and the cost of fast controls $C_T$ verifies
$$C_T\lesssim e^{\frac{K}{\sqrt{T}}}$$
for all $K>3^{1/4}L^{3/2}$.
Moreover, the power of $1/T$ involved in the exponential is optimal.
\end{theorem}

\subsection{Linear KdV equations controlled on the boundary: the case of Dirichlet boundary conditions with a boundary control on the derivative of the state}
In this section, we consider the following controlled linearized KdV equation posed on $(0,T)\times(0,L)$:

\begin{equation}\label{KdVsa}
\begin{aligned}
y_t+y_x+y_{xxx}&=0&\mbox{ in }(0,T)\times(0,L),
\\y(t,0)&=0&\mbox{ in }(0,T),
\\y(t,L)&=0&\mbox{ in }(0,T),
\\y_x(t,L)&=u(t)+y_x(t,0)&\mbox{ in }(0,L),
\end{aligned}
\end{equation}
with initial condition $y^0\in H:=H^{-1}(0,L)$ and control $u\in L^2(0,T)$. We call $A$ the operator $\partial^3_{xxx}+\partial_x$ with domain 
$$\mathcal D(A):=\{y\in H^2(0,L)|y(0)=y(L)=0, y'(0)=y'(L)\}.$$ The eigenvalues are denoted $(i\lambda_n)_{n\in \mathbb Z}$ with $\lambda_n\in\mathbb R$.

This equation describes the propagation of water waves in a uniform channel where $(x,y)$ represents the horizontal and vertical coordinates of the level of water (see for example \cite{MR718811}). We know (see \cite{MR2529319}) that in this case there exists a unique mild solution  $y\in C^0([0,T],H^{-1}(0,L))$. Moreover, it is proved in \cite{MR2529319} that this equation is exactly controllable (and then null controllable) at all time as soon as $L\not\in\mathcal N$ where 
$$\mathcal N:=\{2\pi\sqrt{\frac{k^2+kl+l^2}{3}}|(k,l)\in(\mathbb N^*)^2\}.$$
The original system studied in \cite{97Ro} was 
\begin{equation}\label{KdVnonsa}
\begin{aligned}
y_t+y_x+y_{xxx}&=0&\mbox{ in }(0,T)\times(0,L),
\\y(t,0)&=0&\mbox{ in }(0,T),
\\y(t,L)&=0&\mbox{ in }(0,T),
\\y_x(t,L)&=h(t)&\mbox{ in }(0,L),
\end{aligned}
\end{equation}
with initial condition $y^0\in L^2(0,L)$ and control $h\in L^2(0,T)$.

However, the problem is that the steady-state operator associated to \eqref{KdVnonsa} with the given boundary condition is neither self-adjoint nor skew-adjoint, so we cannot apply directly the results presented before. That is why we have to change a little bit the boundary condition so that the associated steady-state operator becomes skew-adjoint by using the system \eqref{KdVsa} studied in \cite{MR2529319}. 

To be able to apply Theorem~\ref{MTh} or Theorem~\ref{MTh-disp}, we have to study the sequence of eigenvalues $(\lambda_n)_{n\geqslant 1}$. One has the following result:
\begin{lemma}\label{eigenok} $(\lambda_n)_{n\in\mathbb Z}$ is regular and one has
\begin{gather}\label{vpkdv}\lambda_n=\frac{8\pi^3n^3}{L^3}+O(n^2)\end{gather} as $n\rightarrow \pm\infty$.
\end{lemma}
\textbf{Proof of Lemma~\ref{eigenok}.}
This is an immediate consequence of \cite[Proposition 1]{MR2529319}, which gives exactly \eqref{vpkdv} and proves that each eigenspace is of dimension $1$, which implies the regularity of $(\lambda_k)_{k\in\mathbb Z}$ because of the asymptotic behavior given by \eqref{vpkdv}.
\cqfd

From now on, we call $e_k$ one of the unitary eigenvector (for the $H^{-1}$-norm) corresponding to the eigenvalue $i\lambda_k$. We fix an initial condition $y^0:=\sum_{k\in\mathbb Z}a_ke_k\in H^{-1}(0,L)$. As in the previous Subsection, we have for every $\varphi\in \mathcal D(A)$,
$$b(\varphi)=-(\Delta^{-1}\varphi)'(0),$$
so that
$$b=\delta_L'\circ \Delta^{-1},$$
and $$|b_k|=|e_k'(L)|/k^2.$$
To apply Theorem~\ref{MTh}, 
we just need to ensure that
\begin{lemma}\label{bok}
$$b_k\simeq 1.$$
\end{lemma}
\textbf{Proof of Lemma~\ref{bok}.}
$b_k\not=0$ is a consequence of \cite[Lemma 3.5]{97Ro} (because $L\not\in\mathcal N$) and  \cite[Lemma 3.1]{MR2529319} gives immediately that $|e_k'(0)|$ is equivalent as $k\rightarrow \infty$ to $2\pi\sqrt{3}k^2/L^{3/2}$ (because in Lemma 3.1 of \cite{MR2529319} the eigenvectors are normalized in the $L^2$-norm and here in the $H^{-1}$-norm so the behavior of their norm as $k\rightarrow \infty$  has to be multiplied by $k$), so we finally have $b_k\simeq 1$. 
\cqfd
Applying Theorem\ref{MTh-disp}, we obtain directly the following Theorem:
\begin{theorem}\label{cont-kdv-aux}
Let $L\not\in \mathcal N$. Then equation \eqref{KdVsa} is null controllable and the cost of fast controls $C_T$ verifies
$$C_T\leqslant e^{\frac{K}{\sqrt{T}}}$$
for all $K>3^{1/4}L^{3/2}$.
Moreover, the power of $1/T$ involved in the exponential is optimal.
\end{theorem}
\begin{remark}
Using \cite[Remark 1.3]{MR2571687}, one can also add a term of diffusion $-y_{xx}$ in equation \eqref{KdVsa} and obtain the same upper bound as in Theorem~\ref{cont-kdv-aux}.
\end{remark}

\subsection{Anomalous diffusion equation in one dimension}
Let us first consider the $1-D$ Laplace operator $\Delta$ in the domain $D(\Delta):=H^1_0(0,L)$ with state space $H:=H^{-1}(0,L)$. It is well-known that $-\Delta:D(\Delta)\rightarrow H^{-1}(0,L)$ is a definite positive operator with compact resolvent, the $k$-th eigenvalue is 
$$\lambda_k=\frac{k\pi}{L},$$
one of the corresponding normed (in $H$) is 
$$e_k(x):=\frac{\sqrt{2}(1+k^2)^{1/2} \sin(kx)}{\sqrt{L}}.$$
Thanks to the continuous functional calculus for positive self-adjoint operators, one can define any positive power of $-\Delta$. Let us consider here some $\gamma>1/2$ and let us call $\Delta^\gamma:=-(-\Delta)^\gamma$. The domain of $\Delta^\gamma$, that we denote $H_\gamma$, is the completion of $C^\infty_0(0,L)$ for the norm 
$$||\psi||_\gamma:=(\sum_{k\in \mathbb N^*} \frac{<e_k,\psi>_H}{k^{2\gamma}})^{1/2}.$$
We now consider the following equation on $(0,T)\times(0,L)$: 
\begin{equation}\label{anofree}
\begin{aligned}
y_t&=\Delta^\gamma y &\mbox{ in } (0,T)\times(0,L),\\
y(0,.)&=y^0 &\mbox{ in }(0,L) .
\end{aligned}
\end{equation}
This kind of equation can modelize anomaly fast or slow diffusion (see for example \cite{MR2090004}).

We now consider the following controlled equation on $(0,T)\times(0,L)$, that we write under the abstract form
\begin{equation}\label{ano}
\begin{aligned}
y_t&=\Delta^\gamma y+bu&\mbox{ in } (0,T)\times(0,L),\\
y(0,.)&=y^0&\mbox{ in }(0,L), 
\end{aligned}
\end{equation}
where for every $\varphi\in \mathcal D(A)$,
$$b(\varphi)=-(\Delta^{-1}\varphi)'(0),$$
i.e.
$$b:=\delta_0'\circ\Delta^{-1} \in D((-\Delta)^\gamma)'$$
and $u\in L^2(0,T)$. 
If $\gamma \in \mathbb N^*$, one can observe, using integrations by parts, that $b$ corresponds to a boundary control on the left side on the $\gamma-1$-th derivative of $y$, so that $b$ can be considered as a natural extension of the boundary control in the case of non-entire $\gamma$ (this kind of controls has already been introduced in \cite[Section 3.3]{MR2272076} to give some negative results about the control of fractional diffusion equations with $\gamma\leqslant1/2$).

We see that 
$$b_k=|e_k'(L)|/k^2\simeq 1.$$
If $y^0\in H$, then there exists a unique solution of \eqref{ano} in the space $C^0([0,T],H)$ (because $b$ is admissible for the semigroup). To our knowledge, the controllability of anomalous diffusion equations with such a control operator and $\gamma\geqslant 1$ has never been studied before.

Applying directly Theorem~\ref{MTh} and Theorem~\ref{lower-bound}, we obtain:
\begin{theorem}
Assume $\gamma\geqslant 1$.
\label{anoth}
Then Equation \eqref{ano} is null controllable with continuous controls. Moreover, the cost of the control in $L^\infty$ norm, still denoted $C_T$ here, is such that 
$$C_T\leqslant e^{\frac{K}{T^{1/(2\gamma-1)}}}$$ for all $$K>\frac{3(2)^{1/(2\gamma-1)}\pi^{2\gamma/(2\gamma-1)}}{((2\sin(\pi/(4\gamma)))^{2\gamma/(2\gamma-1)})}.$$
Moreover, the power of $1/T$ involved in the exponential is optimal.
\end{theorem}
\subsection{Fractional Schrödinger equation in one dimension}
We keep the notations of the previous subsection.
Let us consider the following fractional Schrödinger equation defined on $(0,T)\times(0,L)$ controlled on one side:
\begin{equation}\label{anosc}
\begin{aligned}
y_t=i\Delta^\gamma y+bu&\mbox{ in } (0,T)\times(0,L),\\
y(0,.)=y^0&\mbox{ in } (0,L),
\end{aligned}
\end{equation}

with initial condition $y^0\in H$, $\gamma>3/4$ and (as in the previous subsection) $b:=\delta_0'\circ\Delta^{-1}$. The $1-D$ Laplace operator $\partial_{xx}$ is considered in the domain $D(\partial_{xx}):= H^1_0(0,L)$. Equation \eqref{ano} has a unique solution in $C^0([0,T],H)$ with $H=H^{-1}(0,L)$. This equation has a physical meaning and can be used to study the energy spectrum of a $1-D$ fractional oscillator or for some fractional Bohr atoms, see for \cite{MR1755089},\cite{MR1948569} or \cite{MR2258580}. As far as we know, the control of this kind of equations has never been studied.
As before,
$|b_k|=|e_k'(0)|/k^2\simeq 1$.
Applying directly Theorem~\ref{MTh}, we obtain:
\begin{theorem}
Assume $\gamma\geqslant 1$.
\label{anoscth}
Then Equation \eqref{anosc} is null controllable. Moreover, the cost of the control $C_T$ is such that 
\begin{gather*}
C_T\lesssim e^{\frac{K}{T^{1/(2\gamma-1)}}}\mbox{  for all }K>\frac{3(2)^{1/(2\gamma-1)}L^{2\gamma/(2\gamma-1)}}{4((\sin(\pi/(2\gamma)))^{2\gamma/(2\gamma-1)}}.
\end{gather*}
Moreover, the power of $1/T$ involved in the exponential is optimal.
\end{theorem}
\subsection*{Acknowledgments}
The author would like to thank Professor Jean-Michel Coron for having attracted his interest on this problem and Ivonne Rivas for corrections and remarks that improved significantly the paper.
\bibliographystyle{plain}
\bibliography{Biblio}

\def\cprime{$'$} \newcommand{\SortNoop}[1]{}
\begin{thebibliography}{10}

\bibitem{MR2851682}
Farid Ammar-Khodja, Assia Benabdallah, Manuel Gonz{\'a}lez-Burgos, and Luz
  de~Teresa.
\newblock The {K}alman condition for the boundary controllability of coupled
  parabolic systems. {B}ounds on biorthogonal families to complex matrix
  exponentials.
\newblock {\em J. Math. Pures Appl. (9)}, 96(6):555--590, 2011.

\bibitem{MR718811}
Jerry Bona and Ragnar Winther.
\newblock The {K}orteweg-de {V}ries equation, posed in a quarter-plane.
\newblock {\em SIAM J. Math. Anal.}, 14(6):1056--1106, 1983.

\bibitem{MR2529319}
Eduardo Cerpa and Emmanuelle Cr{\'e}peau.
\newblock Rapid exponential stabilization for a linear {K}orteweg-de {V}ries
  equation.
\newblock {\em Discrete Contin. Dyn. Syst. Ser. B}, 11(3):655--668, 2009.

\bibitem{MR2302744}
Jean-Michel Coron.
\newblock {\em Control and nonlinearity}, volume 136 of {\em Mathematical
  Surveys and Monographs}.
\newblock American Mathematical Society, Providence, RI, 2007.

\bibitem{05aa}
Jean-Michel Coron and Sergio Guerrero.
\newblock Singular optimal control: a linear 1-{D} parabolic-hyperbolic
  example.
\newblock {\em Asymptot. Anal.}, 44(3-4):237--257, 2005.

\bibitem{71FR}
Hector~O. Fattorini and David~L. Russell.
\newblock Exact controllability theorems for linear parabolic equations in one
  space dimension.
\newblock {\em Arch. Rational Mech. Anal.}, 43:272--292, 1971.

\bibitem{2010-Glass-JFA}
Olivier Glass.
\newblock A complex-analytic approach to the problem of uniform controllability
  of a transport equation in the vanishing viscosity limit.
\newblock {\em J. Funct. Anal.}, 258(3):852--868, 2010.

\bibitem{MR2463799}
Olivier Glass and Sergio Guerrero.
\newblock Some exact controllability results for the linear {K}d{V} equation
  and uniform controllability in the zero-dispersion limit.
\newblock {\em Asymptot. Anal.}, 60(1-2):61--100, 2008.

\bibitem{MR2571687}
Olivier Glass and Sergio Guerrero.
\newblock Uniform controllability of a transport equation in zero
  diffusion-dispersion limit.
\newblock {\em Math. Models Methods Appl. Sci.}, 19(9):1567--1601, 2009.

\bibitem{MR2724598}
Olivier Glass and Sergio Guerrero.
\newblock Controllability of the {K}orteweg-de {V}ries equation from the right
  {D}irichlet boundary condition.
\newblock {\em Systems Control Lett.}, 59(7):390--395, 2010.

\bibitem{MR805273}
Edgardo~N. G{\"u}ichal.
\newblock A lower bound of the norm of the control operator for the heat
  equation.
\newblock {\em J. Math. Anal. Appl.}, 110(2):519--527, 1985.

\bibitem{MR2258580}
Xiaoyi Guo and Mingyu Xu.
\newblock Some physical applications of fractional {S}chr\"odinger equation.
\newblock {\em J. Math. Phys.}, 47(8):082104, 9, 2006.

\bibitem{MR704478}
Lop~F. Ho and David~L. Russell.
\newblock Admissible input elements for systems in {H}ilbert space and a
  {C}arleson measure criterion.
\newblock {\em SIAM J. Control Optim.}, 21(4):614--640, 1983.

\bibitem{MR1755089}
Nikolai Laskin.
\newblock Fractional quantum mechanics and {L}\'evy path integrals.
\newblock {\em Phys. Lett. A}, 268(4-6):298--305, 2000.

\bibitem{MR1948569}
Nikolai Laskin.
\newblock Fractional {S}chr\"odinger equation.
\newblock {\em Phys. Rev. E (3)}, 66(5):056108, 7, 2002.

\bibitem{MR2956149}
Pierre Lissy.
\newblock A link between the cost of fast controls for the 1-{D} heat equation
  and the uniform controllability of a 1-{D} transport-diffusion equation.
\newblock {\em C. R. Math. Acad. Sci. Paris}, 350(11-12):591--595, 2012.

\bibitem{MR2090004}
Ralf Metzler and Joseph Klafter.
\newblock The restaurant at the end of the random walk: recent developments in
  the description of anomalous transport by fractional dynamics.
\newblock {\em J. Phys. A}, 37(31):R161--R208, 2004.

\bibitem{MR2248170}
Sorin Micu and Enrique Zuazua.
\newblock On the controllability of a fractional order parabolic equation.
\newblock {\em SIAM J. Control Optim.}, 44(6):1950--1972 (electronic), 2006.

\bibitem{2004-Miller}
Luc Miller.
\newblock Geometric bounds on the growth rate of null-controllability cost for
  the heat equation in small time.
\newblock {\em J. Differential Equations}, 204(1):202--226, 2004.

\bibitem{MR2062431}
Luc Miller.
\newblock How violent are fast controls for {S}chr\"odinger and plate
  vibrations?
\newblock {\em Arch. Ration. Mech. Anal.}, 172(3):429--456, 2004.

\bibitem{MR2272076}
Luc Miller.
\newblock On the controllability of anomalous diffusions generated by the
  fractional {L}aplacian.
\newblock {\em Math. Control Signals Systems}, 18(3):260--271, 2006.

\bibitem{MR2679651}
Luc Miller.
\newblock A direct {L}ebeau-{R}obbiano strategy for the observability of
  heat-like semigroups.
\newblock {\em Discrete Contin. Dyn. Syst. Ser. B}, 14(4):1465--1485, 2010.

\bibitem{MR2723248}
Frank W.~J. Olver, Daniel~W. Lozier, Ronald~F. Boisvert, and Charles~W. Clark,
  editors.
\newblock {\em N{IST} handbook of mathematical functions}.
\newblock U.S. Department of Commerce National Institute of Standards and
  Technology, Washington, DC, 2010.
\newblock With 1 CD-ROM (Windows, Macintosh and UNIX).

\bibitem{97Ro}
Lionel Rosier.
\newblock Exact boundary controllability for the {K}orteweg-de {V}ries equation
  on a bounded domain.
\newblock {\em ESAIM Control Optim. Calc. Var.}, 2:33--55 (electronic), 1997.

\bibitem{MR2084328}
Lionel Rosier.
\newblock Control of the surface of a fluid by a wavemaker.
\newblock {\em ESAIM Control Optim. Calc. Var.}, 10(3):346--380 (electronic),
  2004.

\bibitem{MR1214759}
David~L. Russell and Bing~Yu Zhang.
\newblock Controllability and stabilizability of the third-order linear
  dispersion equation on a periodic domain.
\newblock {\em SIAM J. Control Optim.}, 31(3):659--676, 1993.

\bibitem{MR1360229}
David~L. Russell and Bing~Yu Zhang.
\newblock Exact controllability and stabilizability of the {K}orteweg-de
  {V}ries equation.
\newblock {\em Trans. Amer. Math. Soc.}, 348(9):3643--3672, 1996.

\bibitem{MR0015553}
Laurent Schwartz.
\newblock Approximation d'un fonction quelconque par des sommes
  d'exponentielles imaginaires.
\newblock {\em Ann. Fac. Sci. Univ. Toulouse (4)}, 6:111--176, 1943.

\bibitem{MR0014502}
Laurent Schwartz.
\newblock {\em \'{E}tude des sommes d'exponentielles r\'eelles}.
\newblock Actualit\'es Sci. Ind., no. 959. Hermann et Cie., Paris, 1943.

\bibitem{MR743923}
Thomas~I. Seidman.
\newblock Two results on exact boundary control of parabolic equations.
\newblock {\em Appl. Math. Optim.}, 11(2):145--152, 1984.

\bibitem{TT}
G\'erald Tenenbaum and Marius Tucsnak.
\newblock New blow-up rates for fast controls of {S}chr\"odinger and heat
  equations.
\newblock {\em J. Differential Equations}, 243(1):70--100, 2007.

\bibitem{MR2859866}
G{\'e}rald Tenenbaum and Marius Tucsnak.
\newblock On the null-controllability of diffusion equations.
\newblock {\em ESAIM Control Optim. Calc. Var.}, 17(4):1088--1100, 2011.

\bibitem{MR920808}
George Weiss.
\newblock Admissibility of input elements for diagonal semigroups on {$l^2$}.
\newblock {\em Systems Control Lett.}, 10(1):79--82, 1988.

\end{thebibliography}
\end{document}